\documentclass[titlepage,12pt]{article}
\usepackage{amssymb}
\usepackage{amsfonts}
\begin{document}


{\Large \bf Beyond Topologies, Part I} \\

{\bf Elem\'{e}r E Rosinger, Jan Harm van der Walt} \\
Department of Mathematics \\
and Applied Mathematics \\
University of Pretoria \\
Pretoria \\
0002 South Africa \\
eerosinger@hotmail.com \\

{\bf Abstract} \\

Arguments on the need, and usefulness, of going beyond the usual Hausdorff-Kuratowski-Bourbaki,
or in short, HKB concept of topology are presented. The motivation comes, among others, from
well known {\it topological type processes}, or in short TTP-s, in the theories of Measure,
Integration and Ordered Spaces. These TTP-s, as shown by the classical characterization given
by the {\it four Moore-Smith conditions}, can {\it no longer} be incorporated within the usual
HKB topologies. One of the most successful recent ways to go beyond HKB topologies is that
developed in Beattie \& Butzmann. It is shown in this work how that extended concept of
topology is a {\it particular} case of the earlier one suggested and used by the first author
in the study of generalized solutions of large classes of nonlinear partial differential
equations. \\ \\

{\large \bf 1. Introduction} \\

{\bf Some of the Main Facts} \\

The starting observation of the approach to pseudo-topologies presented here and introduced
earlier in Rosinger [1-7] is the difference between {\it rigid}, and on the other hand, {\it
nonrigid} mathematical structures, a difference which is explained in the sequel. \\

As it happens, the usual concept of HKB topology is a {\it rigid} mathematical structure, and
this is precisely one of the main reasons for a number of its important deficiencies. In
contradistinction to that, the pseudo-topologies we deal with are {\it nonrigid} mathematical
structures, and thus turn out to have a convenient flexibility. \\

There are two novelties in the way the usual HKB concept of topology is extended in this
work. \\

First, the usual sequences, nets or filters used in describing limits, convergence, etc., are
replaced with pre-ordered sets. This extension proves to be quite natural, and not excessive,
since under rather general conditions, it can - if so required - lead back to the usual
sequences, nets or filters. \\
Second, the basic concept is {\it not} that of convergence, but one of {\it Cauchy-ness} that
models in a rather extended fashion the property of {\it Cauchy} sequences, nets or filters.
Of course, convergence is recovered as a particular case of Cauchy-ness, and the advantage of
such an approach is that {\it completion} becomes easily available in the general case. \\

To be more precise, let us consider the case when filters are used. Then convergence on a
space $X$ is typically a relationship between some filters ${\cal F}$ on $X$ and corresponding
points $x$ in $X$. \\
On the other hand, the mentioned Cauchy-ness is a {\it binary relation} $\Xi$ between certain
pairs of filters on $X$. And a filter ${\cal F}$ on $X$ is Cauchy, if it relates to {\it
itself} according to that binary relation, that is, if ${\cal F} ~\Xi~ {\cal F}$. \\
The particular case of a filter ${\cal F}$ convergent to a point $x$ is recovered when that
filter is in the relation $\Xi$ with the ultrafilter generated by $x$, or more generally, with
some filters closely related to that ultrafilter. \\

The point however is that instead of filters, sequences or nets, one uses this time
pre-ordered sets. Furthermore, instead of starting with convergence, one starts with the
concept of Cauchy-ness. \\

As it turns out, the extension of the usual HKB concept of topology presented here is such
that it can {\it no longer} be contained within the usual Eilenberg - Mac Lane concept of
category, Rosinger [28,29]. In other words, the totality of pseudo-topologies dealt with here
no longer constitutes a usual Eilenberg - Mac Lane category, but one in a more general sense.
The reason for that, however, is not in the possibly excessive nature of the extension of the
HKB concept of topology implemented here, but it is more simple. Namely, it is due to the fact
that the HKB concept of topology is rigid, while the concept of pseudo-topology dealt with
here is nonrigid. \\
This issue, however, will not be pursued here, and respective details can be found in Rosinger
[29]. \\

{\bf Pseudo-Topologies are Nonrigid Structures} \\

Let us note that, as pointed out in Rosinger [2, p. 225], in the case of the usual HKB concept
of topology all the topological entities, such as open sets, closed sets, compact sets,
convergence, continuous functions, etc., are {\it uniquely} determined by {\it one single}
such entity, for instance, that of open sets. In other words, if one for example chooses to
start with the open sets, then all the other topological entities can be defined in a unique
manner based on the given open sets. \\

In contradistinction to that usual situation encountered with the HKB topologies, and as seen
in the sequel, in pseudo-topological structures there is a {\it relative independence} between
various topological entities. Namely

\begin{itemize}

\item the connection between topological entities is no longer rigid to the extreme as in the
case of the usual HKB concept of topology, where one of the topological entities determines in
a unique manner all the other ones; instead, the various topological entities are only
required to satisfy certain compatibility conditions, Rosinger [2, p. 225],

\item there exist, however, several topological entities of prime importance, among them, that
of {\it Cauchy-ness} which is a binary relation between two arbitrary sequences, nets, filters,
etc., and it was first introduced in Rosinger [2], see also Rosinger [1,3-7]

\end{itemize}

In this way, the usual HKB concept of topology can be seen as a {\it rigid} mathematical
structure, while the pseudo-topologies dealt with in this work are {\it nonrigid}. \\

Needless to say, there are many other nonrigid mathematical structures, such a rings, spaces
with measure, topological groups, topological vector spaces, and so on, Rosinger [29, pp.
8,9]. \\

Obviously, an important advantage of a rigid mathematical structure, and in particular, of the
usual HKB concept of topology, is a simplicity of the respective theoretical development. Such
simplicity comes from the fact that one can start with only one single entity, like for
instance the open sets in the case of HKB topologies, and then based on that first entity, all
the other entities can be defined or constructed in a clear and unique manner. \\
Consequently, the impression may be created that one has managed to develop a kind of
universal theory, universal in the sense that there may not be any need for alternative
theories in the respective discipline, as for instance is often the perception about the HKB
topology. \\

The disadvantage of a rigid mathematical structure is in a consequent built in lack of
flexibility regarding the interdependence of the various entities involved, since each of them,
except for a single starting one, are determined uniquely in terms of that latter. And in the
case of the HKB topologies this is manifested, among others, in the difficulties related to
dealing with suitable topologies on spaces of continuous functions, as seen in the sequel. \\

Nonrigid mathematical structures, and in particular, pseudo-topologies, can manifest fewer
difficulties coming from a lack of flexibility. \\

A disadvantage of such nonrigid mathematical structures - as for instance with various
approaches to pseudo-topologies - is in the large variety of ways the respective theories can
be set up. Also, their respective theoretical development may turn out to be more complex than
is the case with rigid mathematical structures. \\
Such facts can lead to the impression that one could not expect to find a universal enough
nonrigid mathematical structure, and for instance, certainly not in the realms of topological
type structures, or in short, TTS-s, and in particular, not in the case of pseudo-
topologies. \\

As it happens so far in the literature on pseudo-topologies, there seems not to be a wider and
explicit enough awareness about the following two facts

\begin{itemize}

\item one should rather use nonrigid structures in order to avoid the difficulties coming from
the lack of flexibility of the rigid concept of usual HKB topology,

\item the likely consequence of using nonrigid structures is the lack of a sufficiently
universal concept of pseudo-topology.

\end{itemize}

As it happens, such a lack of awareness leads to a tendency to develop more and more general
concepts of pseudo-topology, hoping to reach a sufficiently universal one, thus being able to
replace once and for all the usual HKB topology with "THE" one and only "winning" concept of
pseudo-topology. \\
Such an unchecked search for increased generality, however, may easily lead to rather meagre
theories. \\

It also happens in the literature that, even if mainly intuitively, when setting up various
concepts of pseudo-topology a certain restraint is manifested when going away from a rigid
theory towards some nonrigid ones. And certainly, the reason for such a restraint is that one
would like to hold to the advantage of rigid theories which are more simple to develop than
the nonrigid ones. \\

After some decades, the literature on pseudo-topologies appears to have settled in some of its
main trends. And one of them is the preference to start with formalizing in a rather large
variety of ways the concept of {\it convergence}. \\
In this regard it is worth noting that, when back in 1914, Hausdorff created the modern
concept of topology, metric spaces, with their open sets, were taken as a starting point for
generalization. And obviously, a natural way to generalize metric spaces is to leave aside the
metric, but keep the concept of open sets. \\
Unfortunately, the resulting HKB concept of topology proved to have a number of important
deficiencies. And they were manifested not so much on the level of a given topological space,
as on the level of continuous mappings between two arbitrary topological spaces. \\

This points to the fact that the HKB concept of topology tends to fail in a {\it categorial}
sense, that is, through its {\it morphisms}, rather than through its {\it objects}. And by now
this categorial level of failure is quite well understood and formulated by pointing out the
fact that the category of usual HKB topologies is {\it not} Cartesian closed, Herrlich. \\

From this point of view the tendency to base notions of pseudo-topology on concepts of
convergence seems indeed to be more deep and sophisticated than the traditional basing of the
usual HKB concept of topology on open sets. Indeed, unlike open subsets, convergence does
inevitably involve {\it nontrivial} morphisms. Namely, mappings from pre-ordered directed sets
to the respective pseudo-topological spaces. \\

{\bf Limits of Axiomatic Theories} \\

In modern Mathematics it is "axiomatic" that theories are built as {\it axiomatic systems}. \\
Unfortunately however, ever since the early 1930s and G\"{o}del's Incompleteness Theorem, we
cannot disregard the deeply inherent limitations of axiomatic mathematical theories. \\
And that limitation cannot be kept away from nonrigid mathematical structures either, since
such structures are also built as axiomatic theories. \\

Consequently, various theories of pseudo-topologies, including the one dealt with here, are
quite likely bound to suffer from limitations when trying to model in their ultimately
inclusive extent such nontrivial mathematical phenomena like topological type structures, or
in short, TTS-s. \\
After all, G\"{o}del's incompleteness already is manifested in the Peano
axiomatic construction of the natural numbers $\mathbb{N}$. And obviously, TTS-s do not by any
means seem to be simpler than the story of $\mathbb{N}$ ... \\

A possible further development, away from such limitations, has recently appeared, even if it
is not yet clearly in the awareness of those dealing with topology and pseudo-topology. Namely,
{\it self-referential} axiomatic mathematical theories are being developed, Barwise \& Moss.
The novelty - until recently considered to be nothing short of a sheer scandal - in such
theories is in the use of concepts given by definitions which involve a {\it vicious circle},
that is, are self-referential. \\

It may, however, happen that a further better understanding and modelling of TTS-s may benefit
from such a truly novel approach ... \\

{\bf Origins of the Usual Concept of Topology} \\

Our present usual concept of {\it topology} was first formulated in modern terms by Felix
Hausdorff in his 1914 book "Grundz\"{u}ge der Mengenlehre". Hausdorff built on the earlier
work of M Fr\'{e}chet. During the next two decades, major contributions to the establishment
of that abstract or general concept of topology in its present form have been made by LEJ
Brouwer, K Kuratowski, and the Bourbaki group of mathematicians. The Bourbaki group completed
that topological venture with the introduction of the concept of {\it uniform} topology,
towards the end of the 1930s. That was a highly important {\it particular} case of topology,
yet extending fundamental properties of metric spaces, such as for instance, the construction
of the topological {\it completion} of a space. \\

{\bf Difficulties with the Usual Concept of Topology} \\

Strangely enough, soon after, serious {\it limitations} and {\it deficiencies} of that
Hausdorff-Kuratowski-Bourbaki, or in short, HKB concept of topology started to surface. \\

One of the more important such problems came from the shockingly surprising difficulties in
setting up suitable topologies on important and frequently used {\it spaces of functions}. In
this regard most simple {\it linear} situations would already give highly worrying signals. \\

Let for instance $E$ be a locally convex topological space on $\mathbb{R}$, and let $E^*$ be
its topological dual made up of all the continuous linear functionals from $E$ to $\mathbb{R}$.
Let $E^*$ be endowed with any locally convex topology. We consider the usual {\it evaluation}
mapping defined by duality, namely \\

$~~~~~~ ev : E \times E^* \ni ( x, x^* ) \longmapsto\,\, <\, x, x^* \,>\,\,
                                                =~ x^* ( x ) \in \mathbb{R} $ \\

Then rather surprisingly - and also, most inconveniently - it follows that the {\it joint
continuity} of this evaluation mapping $ev$ will actually {\it force} the locally convex
topology on $E$ to be so particular, as to be {\it normable}. \\
Indeed, let us assume that the evaluation mapping $ev : E \times E^* ~\longrightarrow~
\mathbb{R}$ is jointly continuous. Then there exist neighbourhoods $U \subseteq E,~ U^*
\subseteq E^*$ of $0 \in E$ and $0 \in E^*$, respectively, such that $ev ( U, U^* ) \subseteq
[ -1, 1]$. But then $U$ is contained in the polar of $U^*$, therefore $U$ is bounded in $E$.
And since $E$ admits a bounded neighbourhood $U$, it follows that $E$ is normable. \\

There have also been other somewhat less shocking, but on the other hand more pervasive
instances of the deficiencies of the HKB concept of topology. \\

Measure and Integration Theory has been the long ongoing source of some of them. Indeed, about
a decade prior to Hausdorff's introduction of the modern concept of topology, H Lebesgue
established his "dominated convergence" theorem which plays a fundamental role in Lebesgue
integration, and it has no similarly strong version in the case of the Riemann integral. \\
This theorem, involving infinite sequences of integrable functions, their limits, as well as
the limits of their integrals, is obviously about certain {\it topological type processes}, or
in short, TTP-s, in the space of integrable functions. \\
Similar is the situation with Lusin's theorem about the approximation of measurable functions
by continuous ones, or with Egorov's theorem about the almost uniform convergence of
point-wise convergent sequences of measurable functions. \\
Yet none of these three theorems has ever been expressed in terms of HKB topologies. And the
reason - although apparently hardly known well enough, and seldom, if ever, stated explicitly
- is that, simply, they {\it cannot} be expressed in such a manner, due to the fact that they
{\it fail} the Moore-Smith criterion, mentioned in the sequel. \\

Connected with measure theory, one can also mention such a simple and basic concept like
convergence almost everywhere which, however, proves in general {\it not} to be expressible in
terms of HKB topologies, Ordman, since it does not satisfy the mentioned Moore-Smith
criterion. \\

There have been several notable early studies related to extensions of the usual HKB concept
of topology, Choquet, Fischer. However, the respective problems with HKB topologies did not
receive much attention until recent times, Beattie \& Butzmann. \\

Topological type structures, or in short, TTS-s, on ordered spaces have been another rather
considerable source bringing to the fore the limitations of the HKB concept of topology. For
instance, on ordered spaces there is a wealth of natural and useful concepts of convergence
which - in view of the same Moore-Smith criterion - cannot be incorporated within HKB
topologies, Luxemburg \& Zaanen, Zaanen. \\
On the other hand, the particular usefulness of such order based topological structures is
obvious. One of such examples is the 1936 "spectral theorem" of Freudenthal, mentioned later.
A recent remarkable example can be seen in Van der Walt, and for further details in this
regard see also Rosinger [27]. \\

Recently, in connection with attempts in establishing a systematic topological and
differential study of infinite dimensional smooth manifolds, the HKB concept of topology has
met yet another challenging alternative in what came to be called a "convenient setting",
Kriegl \& Michor. That approach - even if so far has not attained all of its main objectives -
can nevertheless quite clearly show the extent of the need to go beyond the HKB concept of
topology. \\

Further relevant facts and arguments related to the need for an extension of the usual HKB
concept of topology are presented in Beattie \& Butzmann [pp. xi-xiii], Beattie, as well as in
the literature mentioned there. \\

{\bf Spaces of Functions and Continuous Convergence} \\

There is a notorious, even if less often mentioned, difficulty with finding suitable HKB
topologies on spaces of functions. And this already happens with linear mappings between
vector spaces. In this regard, the case of the evaluation mapping \\

$~~~~~~ ev : E \times E^* \ni ( x, x^* ) \longmapsto\,\, <\, x, x^* \,>\,\,
                                                =~ x^* ( x ) \in \mathbb{R} $ \\

mentioned above is just one of the simpler, even if critically inconvenient, examples. \\

Here we recall some of the difficulties with the HKB concept of topology related to spaces of
continuous linear mappings, and do so in a more general manner, as pointed out in Rosinger [2,
pp. 223,224]. \\

Let $E$ and $F$ be two locally convex topological vector spaces over $\mathbb{R}$, with their
topology generated respectively by the families of semi-norms $( p_i )_{i \in I}$ and
$( q_j )_{j \in J}$. \\
In such a case it would be most convenient in many situations if the set ${\cal L} ( E, F )$
of linear and continuous mappings from $E$ to $F$ could be obtained in a simple and explicit
manner from the family of sets of linear continuous mappings ${\cal L} ( E_i, F_j )$, with
$i \in I, j \in J$, where $E_i$ and $F_j$ denote respectively the topologies on $E$ and $F$,
generated by the semi-norms $p_i$ and $q_j$. Indeed, the structure of these sets ${\cal L}
( E_i, F_j )$ of linear continuous mappings is perfectly well understood, since they act
between semi-normed spaces. \\

Now, since the topologies on $E$ and $F$ are respectively given by $\sup_{i \in I} E_i$ and
$\sup_{j \in J} F_j$, we obviously have the relations \\

$~~~ \bigcup_{\, i \in I} {\cal L} ( E_i, F ) ~\subseteq~ {\cal L} ( E, F ) ~\subseteq~
                   \bigcap_{\, j \in J} {\cal L} ( E, F_j ) $ \\

$~~~ {\cal L} ( E_i, F ) ~\subseteq~
                   \bigcap_{\, j \in J} {\cal L} ( E_i, F_j ),~~~ i \in I $ \\

$~~~ {\cal L} ( E, F_j ) ~\supseteq~
                   \bigcup_{\, i \in I} {\cal L} ( E_i, F_j ),~~~ j \in J $ \\

Or to argue more simply, let us assume that $J$ has one single element $j$. Then it follows
immediately that \\

$~~~ {\cal L} ( E, F ) ~=~ {\cal L} ( \sup_{i \in I} E_i, F_j ) ~\supseteq~
                                          \bigcup_{i \in I} {\cal L} ( E_i, F_j ) $ \\

Similarly, if $I$ has only one element $i$, then \\

$~~~ {\cal L} ( E, F ) ~=~ {\cal L} ( E_i, \sup_{j \in J} F_j ) ~\subseteq~
                                          \bigcap_{j \in J} {\cal L} ( E_i, F_j ) $ \\

Thus in conclusion, in the general case, the natural relation to expect would be \\

($\cup\cap$)  $~~~ {\cal L} ( E, F ) ~=~
                     \bigcup_{i \in I} \bigcap_{j \in J} {\cal L} ( E_i, F_j ) $ \\

or alternatively \\

($\cap\cup$) $~~~ {\cal L} ( E, F ) ~=~
                       \bigcap_{j \in J} \bigcup_{i \in I} {\cal L} ( E_i, F_j ) $ \\

However, in general, none of these two relations holds, since the connection between ${\cal L}
( E, F )$, and on the other hand, the family ${\cal L} ( E_i, F_j )$, with $i \in I, j \in J$,
turns out typically to be far more involved, thus also far less explicit or simple. \\

This failure to have relations like ($\cup\cap$) or ($\cap\cup$) in general is precisely the
problem address by Marinescu [1,2], although the respective formulation is somewhat different
and less general than the one above. And the suggested solution only regards certain frequent
and useful particular cases, by introducing a suitable concept of {\it pseudo-topology} as an
extension of the HKB concept of topology. \\

One of the remarkable and surprisingly useful concepts in Beattie \& Butzmann is that of the
{\it continuous convergence}. And it addresses precisely the problem of finding suitable
topological type structures on spaces of functions, Beattie \& Butzmann [25-42]. \\
An important consequence of using this concept of continuous convergence is the far more
convenient {\it duality} theory it allows for a large class of vector spaces, Beattie \&
Butzmann [chap. 4], Beattie. \\

{\bf Category Language Helps in Understanding} \\

Put in short terms, a major {\it failure} of the HKB concept of topology is that the
corresponding {\it category} of topological spaces is {\it not} Cartesian closed, Herrlich. In
other words, let $X, Y$ and $Z$ be three spaces with respective usual HKB topologies, and let
endow $X \times Y$ with the corresponding product topology. Set theoretically we have
the relation \\

(EXP) $~~~ Z^{X \times Y} ~=~ ( Z^X )^Y $ \\

which means the existence of the {\it bijective} mapping \\

$~~~~~~ Z^{X \times Y} \ni f ~~\longmapsto~~ f_{ev} \in ( Z^X )^Y $ \\

where \\

$~~~~~~ ( f_{ev} ( y ) ) ( x ) ~=~ f ( x, y ),~~~ x \in X,~ y \in Y $ \\

On the other hand, one {\it cannot} in general find a usual HKB topology on the space ${\cal
C} ( X, Z )$ of continuous functions from $X$ to $Z$, such that by restricting the above
mapping $f \longmapsto f_{ev}$, one would again obtain a bijection, namely \\

$~~~~~~ {\cal C} ( X \times Y, Z ) \ni f ~~\longmapsto~~
                f_{ev} \in {\cal C} ( Y, {\cal C} ( X, Z ) ) $ \\

In other words, within the HKB topologies, we {\it cannot} in general obtain the relation \\

(CONT-EXP) $~~~ {\cal C} ( X \times Y, Z ) ~=~ {\cal C} ( Y, {\cal C} ( X, Z ) ) $ \\

which would be the particular continuous version of (EXP). \\

{\bf The Four Moore-Smith Conditions} \\

The early, 1922 paper of Moore \& Smith introduced the systematic use of {\it nets} in
topology, Kelley [pp. 62-83]. As it happened at the time, the motivation came from a study of
summability. \\

Nets with values in a given set $E$ are natural generalizations of usual sequences $s :
\mathbb{N} \longrightarrow E$. Namely, the index set $\mathbb{N}$ is replaced with an {\it
arbitrary} set $I$, which however is assumed to be endowed with a {\it pre-order} relation
$\leq$ that is at the same time {\it directed}. More precisely, a net in the set $E$ is any
mapping $s : I \longrightarrow E$, where the index set $I$ has a binary relation $\leq$ which
is {\it reflexive} and  {\it transitive}, thus it is a {\it pre-order}, and in addition, it
has the {\it directedness} property that $\forall~ i, j \in I : \exists~ k \in I : i, j \leq
k$. In general, $\leq$ need not be antisymmetric as well, namely, one does not necessarily
ask that $\forall~ i, j \in I : (~ i \leq j,~ j \leq i ~) \Longrightarrow i = j $. \\

One of the basic topological issues arising in connection with nets is the following one,
Kelley [pp. 73,74]. \\

Suppose for a certain set $E$, we are given a class ${\cal S}$ of pairs $( s, x )$, where $s$
are nets in $E$, while $x \in E$. \\

The intended meaning of the binary relationship $( s, x ) \in {\cal S}$ is that the net $s$
{\it converges to} $x$ in the sense of ${\cal S}$, and it is denoted by \\

$~~~~~~ \lim_{\,{\cal S}}~ s ~=~ x $ \\

The topological issue is :

\begin{itemize}

\item Which are the {\it necessary} and {\it sufficient} conditions on the class ${\cal S}$,
so that there exists a usual HKB topology $\tau$ on $E$, for which we have the {\it
equivalence} between the usual convergence in topology, and on the other hand, that given by
the class ${\cal S}$, namely

\end{itemize}

\medskip

(1.1)~~~ $ \begin{array}{l}
               \forall~ x \in E,~ s : I \longrightarrow E ~~\mbox{net in}~~ E : \\ \\
               ~~~~ \lim_{\,\tau}~ s ~=~ x ~~~\Longleftrightarrow~~~
                                                \lim_{\,{\cal S}}~ s ~=~ x
           \end{array} $ \\

\bigskip
The {\it characterization} of this equivalence is given by the following {\it four Moore-Smith
conditions} (1.2) - (1.5) on the class ${\cal S}$, Kelley [p. 74], namely \\

(1.2)~~~ $ \begin{array}{l}
               \forall~ x \in E,~ s : I \longrightarrow E ~~\mbox{net in}~~ E : \\ \\
               ~~~~ (~~ s ( i ) ~=~ x, ~~\mbox{for}~ i \in I ~~)
                                   ~~~\Longrightarrow~~~ \lim_{\,{\cal S}}~ x
            \end{array} $ \\ \\ \\

(1.3)~~~ $ \begin{array}{l}
               \forall~ ( s , x ) \in {\cal S} : \\ \\
               \forall~ t : I \longrightarrow E ~~\mbox{net in}~~ E : \\ \\
               ~~~~ (~~ t ~~\mbox{subnet of}~~ s ~~)
                                     ~~~\Longrightarrow~~~  \lim_{\,{\cal S}}~ t ~=~ x
            \end{array} $ \\ \\ \\

(1.4)~~~ $ \begin{array}{l}
              \forall~ x \in E,~ s : I \longrightarrow E ~~\mbox{net in}~~ E : \\ \\
               ~~~~~~ ( s, x ) \notin {\cal S} ~~~\Longrightarrow~~~
                  \left ( ~~ \begin{array}{l}
                                  \exists~ t ~~\mbox{subnet of}~~ s : \\
                                  \forall~ t^{\,\prime} ~~\mbox{subnet of}~~ t : \\
                                  ~~~~ ( t^{\,\prime}, x ) \notin {\cal S}
                             \end{array} ~~ \right )
             \end{array} $ \\ \\

We note that an equivalent form of condition (1.4) is the following one \\

(1.4$^{\,\prime}$\,)~~~  $ \begin{array}{l}
              \forall~ x \in E,~ s : I \longrightarrow E ~~\mbox{net in}~~ E : \\ \\
                  ~~~~\left ( ~~ \begin{array}{l}
                                  \forall~ t ~~\mbox{subnet of}~~ s : \\
                                  \exists~ t^{\,\prime} ~~\mbox{subnet of}~~ t : \\
                                  ~~~~  \lim_{\,{\cal S}}~ t^{\,\prime} ~=~ x
                             \end{array} ~~ \right ) ~~~\Longrightarrow~~~
                      \lim_{\,{\cal S}}~ s ~=~ x
             \end{array} $ \\ \\

For the formulation of the last, that is, {\it fourth Moore-Smith condition}, we need some
preparation. Let $( \Lambda, \leq )$ be a directed pre-order. Let $( I_\lambda, \leq )$ be a
directed pre-order, for every $\lambda \in \Lambda$. Then on the set \\

$~~~~~~ J ~=~ \Lambda \times \prod_{\, \lambda \in \Lambda}\, I_\lambda $ \\

there is a natural directed pre-order $\leq$ defined by \\

$~~~~~~ ( \lambda, f ) ~\leq~ ( \mu, g ) ~~~\Longleftrightarrow~~~
                    \left( \begin{array}{l}
                             ~~~*)~~ \lambda \leq \mu \\ \\
                             ~**)~~  f ( \nu ) \leq g ( \nu ), ~~\mbox{for}~ \nu \in \Lambda
                             \end{array} ~~\right) $ \\ \\

where $\lambda, \mu \in \Lambda$ and $f, g \in \prod_{\, \nu \in \Lambda}\, I_\nu$. Let us now
define the {\it diagonal} mapping \\

$~~~~~~ \delta : J ~=~ \Lambda \times \prod_{\, \nu \in \Lambda}\, I_\nu
                      \ni ( \lambda, f ) ~\longmapsto~ f ( \lambda ) \in I_\lambda $ \\

Suppose for every $\lambda \in \Lambda$, one is given $( s_\lambda, x_\lambda ) \in {\cal S}$,
where $s_\lambda : I_\lambda \longrightarrow E$ is a net and $x_\lambda \in E$. Then one can
define the net $s : \Lambda \longrightarrow E$ by $s ( \lambda ) = x_\lambda$. Further, one
defines the net $t : J \longrightarrow E$ by $t ( \lambda, f ) = s_\lambda
( \delta ( \lambda, f ) ) = s_\lambda ( f ( \lambda ) )$, for $\lambda \in \Lambda$ and $f \in
\prod_{\, \nu \in \Lambda}\, I_\nu$. \\
Finally, suppose given $x \in E$. \\

Then the {\it fourth Moore-Smith condition} is as follows \\

(1.5)~~~ $ \left(
               \begin{array}{l}
                 ~~~*)~~  \lim_{\,{\cal S}}~ s_\lambda ~=~ x_\lambda,
                                        ~~\mbox{for}~ \lambda \in \Lambda \\ \\
                 ~**)~~ \lim_{\,{\cal S}}~ s ~=~ x
               \end{array} ~~\right) ~~~\Longrightarrow~~~  \lim_{\,{\cal S}}~ t ~=~ x $ \\ \\

A consequence of the above four Moore-Smith conditions is that, in case one or more of them
are {\it not} satisfied by the class ${\cal S}$, then for whichever usual HKB topology on $E$,
there {\it cannot} be precisely the same amount of convergent nets, thus the equivalence (1.1)
does {\it not} hold, since at least one of the two implications "$\Longrightarrow$" or "
$\Longleftarrow$" will fail to be valid. \\

Clearly, it is easy to construct topologies $\tau$ on $E$ such that one or the other of the
implications in (1.1) holds. \\

For instance, we can consider the set ${\cal T}op_{+}$ of all topologies $\tau$ on $E$, for
which we have the implication \\

(1.6)~~~ $ \begin{array}{l}
               \forall~ x \in E,~ s : I \longrightarrow E ~~\mbox{net in}~~ E : \\ \\
               ~~~~ \lim_{\,\tau}~ s ~=~ x ~~~\Longrightarrow~~~  \lim_{\,{\cal S}}~ s ~=~ x
           \end{array} $ \\

Obviously, this set ${\cal T}op_{+}$ is not void, since the finest topology on $E$ belongs to
it. Indeed, in that finest topology, only constant nets converge. And then, according to (1.2),
the implication "$\Longrightarrow$" in (1.6) does hold. \\

Conversely, we can consider the set ${\cal T}op_{-}$ of all topologies $\tau$ on $E$, for
which we have the implication \\

(1.7)~~~ $ \begin{array}{l}
               \forall~ x \in E,~ s : I \longrightarrow E ~~\mbox{net in}~~ E : \\ \\
               ~~~~ \lim_{\,\tau}~ s ~=~ x ~~~\Longleftarrow~~~  \lim_{\,{\cal S}}~ s ~=~ x
           \end{array} $ \\

Then again, this set ${\cal T}op_{-}$ is not void, since the coarsest topology on $E$ belongs
to it. Indeed, in that coarsest topology, all nets converge. And then the implication
"$\Longleftarrow$" in (1.7) trivially holds. \\

It should be mentioned that there are more recent alternative versions of the four Moore-Smith
conditions (1.2) - (1.5), Schechter [pp. 413, 414]. \\

Two typical examples regarding the {\it failure} of the equivalence in (1.1) are given by the
{\it convergence almost everywhere} and the {\it uniform convergence almost everywhere},
Ordman, Cohn, Schechter [p. 564]. \\

For convenience, let us recall here a few details. Let $( X, \Sigma, \mu )$ be a space with a
measure. Let ${\cal M} ( X, \Sigma )$ be the space of real valued measurable functions on $X$.
Further, let $s : I \longrightarrow {\cal M} ( X, \Sigma )$ be any net, and let $f \in {\cal M}
( X, \Sigma )$ be any real valued measurable function. \\

We say that the measurable functions $f_i = s ( i )$, with $i \in I$, {\it converge almost
everywhere} to the measurable function $f \in {\cal M} ( X, \Sigma )$ , if and only if \\

$~~~~~~ \mu^*\, ( X \setminus
    \{ x \in X ~|~ \lim_{\,i \in I}~ f_i ( x ) = f ( x ) \} ) ~=~ 0 $ \\

where $\mu^*$ is the {\it outer measure} generated by $\mu$. \\

Let us now take $E = {\cal M} ( X, \Sigma )$, and define the class ${\cal S}_{ae}$ by the
condition \\

$ ( s, f ) \in {\cal S}_{ae} ~~~\Longleftrightarrow~~~
                  (~~ f_i = s ( i ), ~\mbox{with}~ i \in I,
                         ~~\mbox{converge almost everywhere to}~ f ~~) $ \\

where the measurable functions $f \in E = {\cal M} ( X, \Sigma )$ and the nets $s : I
\longrightarrow E = {\cal M}( X, \Sigma )$ are arbitrary. \\

Then as is well known, Schechter [p. 564], the third of the above Moore-Smith conditions,
namely (1.4), is not satisfied in general. \\

Therefore, there {\it cannot} in general exist any usual HKB topology on the space $E =
{\cal M}( X, \Sigma )$ of measurable functions which would give precisely the same convergence
as the convergence almost everywhere. \\

A similar situation happens, Schechter [p. 564], with uniform convergence almost everywhere,
which is defined as follows. The measurable functions $f_i = s ( i )$, with $i \in I$, {\it
converge uniformly almost everywhere} to the measurable function $f \in {\cal M} ( X,
\Sigma )$ , if and only if \\

$~~~~~~ \begin{array}{l}
           \forall~~ \epsilon > 0 : \\ \\
           \exists~~ X_\epsilon \in \Sigma : \\ \\
           ~~~ f_i ~~\mbox{converges uniformly to}~~ f ~~\mbox{on}~~ X \setminus X_\epsilon
        \end{array} $ \\

A consequence of the above is that the celebrated "dominated convergence" theorem of Lebesgue
{\it cannot} be formulated in terms of HKB topologies. Indeed, this theorem is about sequences
of integrable functions which {\it converge almost everywhere}, and as mentioned above, such a
convergence is beyond the reach of HKB topologies. \\

{\bf Going Beyond the HKB Topology} \\

Nearly a century of parallel developments have taken place in which the usual HKB concept of
topology has been pursued along with a variety of alternative topological type processes and
structures in Measure and Integration Theory, Theory of Ordered Spaces, etc. These parallel
developments, however, have left the issue of a proper extension of the HKB concept of
topology just about as open as it had always been. \\
And it is by now quite obvious that the way forward towards such an extension is not simply by
setting up one or another concept of convergence with possible relaxations of some of the four
Moore-Smith conditions which characterize convergence in usual HKB topologies. \\
In other words, finding a more general concept of topology, and at the same time, not losing
the essential - even if less explicitly manifest - phenomena, aspects, properties, etc.,
involved has proven not to be an easy task. \\

It is, of course, due precisely to this difficulty that a variety of rather different
pseudo-topological attempts have been made in order to go beyond the HKB framework. \\

And as if to highlight the {\it complexity} and {\it depth} of that enterprise, we can recall
that even Category Theory has got involved in it. For instance, as earlier in Herrlich, so
more recently in Clementino et.al., is stated that :

\begin{quote}

"Failure to be Cartesian closed is one of the main defects of the category of topological
spaces."

\end{quote}

\bigskip

{\large \bf 2. Topological Type Structures, or TTS-s} \\

{\bf Towards an Extension of the HKB Concept of Topology} \\

We shall present one of the extensions of the HKB concept of topology which was suggested in
the 1960s, in Rosinger [1-7]. \\
Then we shall show how that extension of the HKB concept of topology incorporates as a {\it
particular} case the more recent - and highly successful - such extended concept of topology
in Beattie \& Butzmann. \\

The motivation for Rosinger [1-7] was given by well known difficulties in Functional Analysis
manifested by the HKB concept of topology. More specifically, these difficulties arose in the
study of generalized solutions of linear and nonlinear partial differential equations. \\
One of the early formulations of these difficulties was in the celebrated 1954 impossibility
result of L Schwartz which was claiming that the space ${\cal D}^{\,\prime}$ of the Schwartz
distributions could not be usefully embedded into differential algebras of generalized
functions. \\

As it turned out not much later, that claim proved to be {\it incorrect}. \\
Indeed, large classes of differential algebras of generalized functions were constructed and
used in the solution of a considerable variety of linear and nonlinear partial differential
equations, as well as in applications to singular differential geometry, Lie groups, and so on,
see Rosinger [5-7, 9-26], Mallios \& Rosinger [1,2], Rosinger \& Walus [1,2], Mallios, and
also the whole subject field :

\smallskip
{~~~~~~~~~~} 46F30 ~~at~~ www.ams.org/msc/46Fxx.html

\smallskip
And all these algebras contain, among other spaces of generalized functions, also the space
${\cal D}^{\,\prime}$ of the Schwartz distributions\\

Over the years, a variety of ways for going beyond the HKB concept of topology have been
suggested in the literature. Some of the earlier ones can be found in  Alfsen \& Njastad,
Choquet, Cs\'{a}sz\'{a}r, Doicinov, Dupont, Fischer, Fr\"{o}licher \& Bucher, Hacque, Haddad,
Hammer, Jarchow [1,2], Leader, Marinescu [1,2], Steiner. More recent references are presented
in Beattie \& Butzmann. \\

One of the ever most seminal and successful extensions of the HKB concept of topology was
presented in Beattie \& Butzmann. \\
As mentioned, however, this extended concept proves to be a {\it particular} case of the {\it
topological type structures}, or in short, TTS-s, introduced and developed in Rosinger [2-7],
and presented briefly in the sequel. \\

{\bf Topological Type Structures and Topological Type Processes} \\

Let us make a further clarification of the terminology. \\

"Topological type structures", or TTS-s, are meant in this work to be the discipline
studied by the usual HKB concept of topology, as well as by its various pseudo-topological
type extensions. And as we have seen, in Measure and Integration Theory, or in the Theory of
Ordered Spaces, there are plenty of important TTS-s which cannot be dealt with conveniently by
the HKB concept of topology. \\

On the other hand, such topological type structures are supposed to be defined in terms of
certain "topological type processes", or TTP-s, such as for instance sequences, nets,
filters, and so on. \\
In this way, topological type processes, or TTP-s, are but the {\it building blocks} of
topological type structures, or TTS-s. \\

{\bf Topological Type Processes, or TTP-s} \\

Let $E$ be any nonvoid set on which we are interested in a {\it topological type structure},
or TTS. \\

As mentioned, a first departure - introduced up-front - from the usual ways encountered in
topology, or for that matter pseudo-topology, appears already in the definition of the {\it
topological type processes}, or TTP-s, on any given set $E$. Namely, instead of considering
for that purpose the usual sequences, nets, filters, etc., on the set $E$, we go one level
higher in abstraction and consider any nonvoid set $E^{\,\prime}$ which may give all the {\it
topological type processes} on $E$. \\
Also as mentioned, this rather general step, however, is tempered by the fact that, without
loss of generality, we can consider such sets $E^{\,\prime}$ of TTP-s as being endowed with a
pre-order. And if and when desired, this pre-order structure proves to be able to lead back to
nets or filters on $E$. \\

The motivation for such a generalization is simple. Indeed, as is well known, usual sequences
of elements in $E$, that is, mappings $s : \mathbb{N} \longrightarrow E$, are in general
insufficient for the description of a large variety of topological type structures on $E$ even
in the particular case of HKB topologies. On the other hand, nets which, as mentioned, are
mappings $s : I \longrightarrow E$, where $I$ is an arbitrary {\it pre-ordered directed} set
of {\it indices}, present technical difficulties. First of all, since they involve arbitrary
sets of indices $I$, the totality of nets on $E$ is no longer a set, but a category. Further,
in case there is an algebraic structure on $E$, such as a group, vector space, algebra, etc.,
it is not easy to consider a similar algebraic structure on the totality of nets on $E$, since
obviously not all such nets have the same set $I$ of indices. Consequently, it is not easy in
terms of arbitrary nets to bring algebra and topology together. \\
Such technical difficulties were in part the reason why filters on $E$ have long been
considered. And such filters have the major advantage to be definable in terms of $E$ alone,
without the need for any other sets, such as for instance, the arbitrary index sets $I$ needed
for nets. As it happens, however, filters on $E$ are again not convenient when there is some
algebraic structure on that set and we want in addition to bring in a topology compatible with
it. Indeed, say, E is a group with the group operation $\circ$. Let $Fil ( E )$ denote the set
of all filters on $E$. Then we can easily and quite naturally extend the operation $\circ$ to
$Fil ( E )$, as follows. Given on $E$ two filters ${\cal F, G} \in Fil ( E )$, we can define
${\cal F} \circ {\cal G} = \{~ F \circ G ~~|~~ F \in {\cal F},~ G \in {\cal G} ~\}$, where as
usual, we defined $F \circ G = \{~ x \circ y ~~|~~ x \in F,~ y \in G ~\}$, taking into account
that $F, G \subseteq E$, in view of the definition of filters on $E$. However, this extended
operation $\circ$ on $Fil ( E )$ will {\it no} longer be a group, unless we are in the trivial
case when $E$ has one single element. \\

In this way, in order to avoid such rather simple but awkward technical difficulties, it
appears to be quite appropriate to consider any suitable, but otherwise quite arbitrary set
$E^{\, \prime}$ as giving all the topological type processes, or TTP-s, which we shall use in
order to define one or another topological type structure, or TTS, on $E$. \\
And as seen later, assuming the existence of a pre-order on $E^{\, \prime}$ is not going to
lead to technical difficulties. \\

{\bf Presence of Natural Pre-Order on Topological Type \\ Processes} \\

It is important to note, Rosinger [7], that in the case of the usual topological type
processes, or TTP-s, namely, when $E^{\,\prime}$ is given by the set of all sequences, nets,
filters, etc., on $E$, there is a natural {\it pre-order} relation $\leq$ on $E^{\,\prime}$,
that is, a binary relation which is {\it reflexive} and {\it transitive}, but need {\it not}
always be as well {\it antisymmetric}. \\
Indeed, sequences have subsequences, nets have subnets, filters have finer filters, and so
on. \\

Therefore, by considering such a pre-order $\leq$ on $E^{\,\prime}$, one can compensate for
the rather high level of abstraction which brings in an arbitrary set $E^{\,\prime}$ as the
topological type processes, or TTP-s, on the space of interest $E$. \\
And as mentioned, such a compensation can be made naturally and always, thus without loss of
generality. \\
We shall deal with this issue in more detail in the sequel. \\

{\bf Cauchy-ness as a the Basic Concept} \\

As second departure, we note that in usual HKB concept of topology the first and most general
concept is {\it not} that of a uniform space, but of a topological one which need not
necessarily be uniform. \\
In other words, the first and most general concept is that of {\it convergent} topological
type processes, TTP-s, be they sequences, nets, filters, etc., and {\it not} of {\it Cauchy}
TTP-s. \\
This specific order of priority in defining concepts creates difficulties when we try to go
beyond the usual HKB concept of topology, since in the case of the corresponding more general
TTS-s defined in terms of convergence it is not so easy to find appropriate extended concepts
of Cauchy-ness. \\

Consequently, we reverse the above usual priority order, and when defining the most general TTS-S, we shall start with
{\it Cauchy} TTP-s, and after that define the {\it convergent} TTP-s, and do so in terms of the already defined Cauchy
TTP-s \\

{\bf Cauchy-ness as a Binary Relation} \\

The third departure is based on a fact which, no matter how obvious at a more careful analysis,
it is nevertheless not always explicitly enough brought to the fore. Let us present the
respective issue in a simple and well known situation. Let $E$ be a space with a metric $d : E
\times E \longrightarrow [ 0, \infty )$. Then, as is well known, a sequence $s : \mathbb{N}
\longrightarrow E$ is called {\it Cauchy}, if and only if \\

~~~~~~$ \begin{array}{l}
            \forall~ \epsilon > 0 : \\
            \exists~ n \in \mathbb{N} : \\
            \forall~ k, l \in \mathbb{N} : \\
            ~~~ k, l \geq n ~~\Longrightarrow~~ d ( x_k, x_l ) \leq \epsilon
         \end{array} $ \\

Now clearly, the property involved, namely \\

~~~~~~$ k, l \geq n ~~\Longrightarrow~~ d ( x_k, x_l ) \leq \epsilon $ \\

is a {\it binary} relation with respect to the terms of the given sequence $s$. And it is
expressed exclusively with reference to the terms of the respective sequence $s$. That is, no
limit of the sequence is involved, even if such a limit may happen to exist. \\

This is in obvious contradistinction with what happens in the case the sequence $s$ is {\it
convergent} in the usual sense. Indeed, $s$ is convergent in the given metric space to some
$x \in E$, if and only if \\

~~~~~~$ \begin{array}{l}
            \forall~ \epsilon > 0 : \\
            \exists~ n \in \mathbb{N} : \\
            \forall~ m \in \mathbb{N} : \\
            ~~~ m \geq n ~~\Longrightarrow~~ d ( x, x_m ) \leq \epsilon
         \end{array} $ \\

And here, the property involved, namely \\

~~~~~~$ m \geq n ~~\Longrightarrow~~ d ( x, x_m ) \leq \epsilon $ \\

is clearly a {\it unary} and {\it not} binary relation with respect to the terms of the
sequence $s$ under consideration. Furthermore, it is {\it not} formulated exclusively in terms of the sequence $s$ alone,
since it refers to its limit $x$ as well. \\

In view of the above, and in view of the fact that here in this work the first basic concept
is Cauchy-ness, and not convergence, we shall define a TTS on the set$E$ by a {\it binary}
relation $\Xi$ on the TTP-s given by $E^{\, \prime}$, namely \\

(2.1)~~~ $ \Xi ~\subseteq~ E^{\, \prime} \times E^{\, \prime} $ \\

The presence of this binary relation $\Xi$ ~- supposed to model the {\it Cauchy-ness} of the
TTP-s in $E^{\, \prime}$ - is in fact the main {\it novelty} of the approach in Rosinger
[2-7]. \\
Given two arbitrary TTP-s $x^{\, \prime}, y^{\, \prime} \in E^{\, \prime}$, the meaning of the
fact that they are in the binary relation $\Xi$, that is \\

(2.2)~~~ $ ( x^{\, \prime}, y^{\, \prime} ) ~\in~ \Xi $ \\

or written equivalently and more simply \\

(2.2$^{\, \prime}$)~~~ $ x^{\, \prime} ~\Xi~ y^{\, \prime} $ \\

is that, in case the space $E$ had a completion $E^{\#}$, then {\it both}~ $x^{\, \prime}$ and
$y^{\, \prime}$ would {\it converge} in that completion to the {\it same} element $x^{\#}$,
see (2.14) in the sequel. \\

Whenever for two TTP-s\, $x^{\, \prime}, y^{\, \prime} \in E^{\, \prime}$ we have the
relationship $x^{\, \prime} ~\Xi~ y^{\, \prime}$, we shall say that $x^{\, \prime}$ and
$y^{\, \prime}$ are {\it Cauchy related}. \\

It is important to note that, in general, the binary relation $\Xi$ is {\it not} supposed to
be reflexive or transitive. However, as seen in (TTS2) below, it is always required to be {\it
symmetric}. Also, when restricted to the Cauchy TTP-s, it is reflexive as well. \\

{\bf Localization of Topological Type Processes} \\

As is well known, in topology it is important to consider sequences, nets, filters, etc.,
which are defined not only on the whole space $E$ under consideration, but also on various of
its nonvoid subsets $A \subseteq E$. \\

The generalization we shall use in this work is given by mappings \\

(2.3)~~~ $ T : {\cal P} ( E ) ~\longrightarrow~ {\cal P} ( E^{\, \prime} ) $ \\

from the power set of the space $E$ under consideration to the power set of $E^{\, \prime}$ of
associated TTP-s. For $A \subseteq E$, the corresponding subset $T ( A ) \subseteq
E^{\, \prime}$ is supposed to represent all the TTP-s based in $A$. For instance, if
$E^{\, \prime}$ is given by all the usual sequences $s : \mathbb{N} \longrightarrow E$, then
$T ( A )$ is the set of all usual sequences with values in $A$, that is, the set of all the
sequences $s : \mathbb{N} \longrightarrow A$. Of course, upon convenience, one may as well
define $T ( A )$ as the larger set \\

~~~~~~$ T ( A ) ~=~ \left \{~ s : \mathbb{N} \longrightarrow E ~~~
                             \begin{array}{|l}
                                       ~~~ \exists~ n \in \mathbb{N} : \\
                                       ~~~ \forall~ m \in \mathbb{N} : \\
                                       ~~~~~ m \geq n ~\Longrightarrow~ s ( m ) \in A
                             \end{array} ~\right \} $ \\

Obviously, if one uses nets, filters, etc., for TTP-s, one can define the mapping $T$
similarly. \\

The {\it two} conditions such a mapping $T$ is supposed to satisfy are \\

(TTP1)~~~ $ \phi \neq A \subseteq B \subseteq E ~~\Longrightarrow~~
                      \phi \neq T ( A ) \subseteq T ( B ) \subseteq E^{\, \prime} $ \\

(TTP2)~~~ $ T ( E ) ~=~ E^{\, \prime} $ \\

The meaning and necessity of condition (TTP1) is obvious. The condition (TTP2) simply means
that in the set $E^{\, \prime}$ of TTP-s associated with $E$ there are {\it no} redundant
elements. In other words, all TTP-s in $E^{\, \prime}$ are based in $E$ since they belong to
$T ( E )$. \\

{\bf Definition of Topological Type Structures} \\

With the above preparation and motivation, we can come now to the general definition of {\it
topological type structures}, or in short, TTS-s,  Rosinger [2-7], which are given by any
{\it quadruplet} \\

(TTS)~~~ $ \sigma ~=~ ( E, E^{\, \prime}, T , \Xi ) $ \\

which satisfies the conditions (TTP1), (TTP2) with respect to the mapping $T$, as well as the
following three additional ones, with respect to the binary relation $\Xi$ of Cauchy-ness,
namely \\

(TTS1)~~~ $ \begin{array}{l}
                    \forall~ x \in E : \\ \\
                    \forall~ x^{\, \prime}, y^{\, \prime} \in T ( \{ x \} ) \subseteq
                                                                     E^{\, \prime} : \\ \\
                    ~~~ x^{\, \prime} ~\Xi~ y^{\, \prime}
             \end{array} $ \\

also \\

(TTS2)~~~ $ \begin{array}{l}
                    \forall~ x^{\, \prime}, y^{\, \prime} \in E^{\, \prime} : \\ \\
                    ~~~ x^{\, \prime} ~\Xi~ y^{\, \prime} ~~\Longrightarrow~~
                                                    y^{\, \prime} ~\Xi~ x^{\, \prime}
            \end{array} $ \\

as well as \\

(TTS3)~~~ $ \begin{array}{l}
                    \forall~ x^{\, \prime}, y^{\, \prime} \in E^{\, \prime} : \\ \\
                    ~~~ x^{\, \prime} ~\Xi~ y^{\, \prime} ~~\Longrightarrow~~
                                                    x^{\, \prime} ~\Xi~ x^{\, \prime}
            \end{array} $ \\

The meaning of these three conditions is as follows. Condition (TTS1) is in fact a rather {\it
trivial} requirement. Indeed, it says that if two TTP-s $x^{\, \prime}, y^{\, \prime} \in
E^{\,\prime}$ are based in the same point $x$ of the space $E$, then they are Cauchy related.
Here we note that TTP-s based in a single point $x \in E$ are usually constant sequences or
nets with the terms equal to $x$, filters generated by the single point $x$, thus fixed
ultrafilters, etc. Therefore, in usual HKB topologies they are certainly convergent to such an
$x$, while in uniform topologies they are also Cauchy. \\
The condition (TTS2) is {\it natural}, since is requires that the binary relation $\Xi$ be
symmetric. And certainly, there seems to exist no a priori reason why one should not assume
that. \\
Finally, the condition (TTS3) is the other {\it novelty}. And it requires that in case two
TTP-s $ x^{\, \prime}, y^{\, \prime} \in E^{\, \prime}$ are Cauchy related to one another,
then each of them is Cauchy related to {\it itself}. The meaning of this condition is further
clarified in the definition of {\it Cauchy TTP-s} given in (2.4) below. \\

{\bf Cauchy and Convergent Topological Type Processes} \\

Given $\sigma ~=~ ( E, E^{\, \prime}, T , \Xi )$ any TTS on $E$, we define the set of {\it
Cauchy TTP-s}~ in $\sigma$ as given by \\

(2.4)~~~ $ Cauchy~ ( \sigma ) ~=~ \{~ x^{\, \prime} \in E^{\, \prime} ~~|~~
                                             x^{\, \prime} ~\Xi~ x^{\, \prime} ~\} $ \\

Further, for any $x \in E$, we define the set of {\it TTP-s convergent to} $x$ in $\sigma$ as
being given by \\

(2.5)~~~ $ Conver~ ( \sigma, x ) ~=~
             \left \{~ x^{\, \prime} \in E^{\, \prime} ~~
                    \begin{array}{|l} ~~
                       \exists~ x^{\, \prime}_0 \in T ( \{ x \} ) \subseteq E^{\, \prime} : \\ \\
                        ~~~~~ x^{\, \prime} ~\Xi~ x^{\, \prime}_0
                    \end{array} ~\right \} $ \\

Clearly, in view of (TTS1), (TTS3), we have the inclusions \\

(2.6)~~~ $ \begin{array}{l}
                   \forall~ x \in E : \\ \\
                   ~~~~ \phi \neq T ( \{ x \} ) ~\subseteq~ Conver~ ( \sigma, x )
                                      ~\subseteq~ Cauchy~ ( \sigma ) ~\subseteq~ E^{\, \prime}
            \end{array} $ \\

which extend the customary relationships between constant, convergent and Cauchy sequences,
nets, filters, etc., in case of usual uniform spaces. \\

We note, see (2.4), that in a TTS $\sigma = ( E, E^{\, \prime}, T , \Xi )$, the binary
relation $\Xi$ is {\it reflexive} precisely on the subset $Cauchy~ ( \sigma )$ of all the
TTP-s $E^{\,\prime}$. \\

{\bf Complete Spaces} \\

From the above follows naturally the definition of complete TTS-s. Namely, given any $\sigma
~=~ ( E, E^{\, \prime}, T , \Xi )$, we say that $\sigma$ is {\it complete}, if and only if \\

(2.7)~~~ $  Cauchy~ ( \sigma ) ~=~ \bigcup_{x \in E}~ Conver~ ( \sigma, x ) $ \\

In view of (2.6), the inclusion $"\supseteq"$ in (2.7) always holds. \\

{\bf Topological Support} \\

Given $\sigma ~=~ ( E, E^{\, \prime}, T , \Xi )$ any TTS on $E$, it will often be convenient
to consider its {\it topological support on} $E$ \\

(2.8)~~~ $ ( E, E^{\, \prime}, T  ) $ \\

Clearly, such a topological support on a set $E$ can be defined or considered independently of
the binary relation $\Xi \subseteq E^{\,\prime} \times E^{\,\prime}$. \\
Also, a given topological support (2.8) can be associated with a variety of binary relations
$\Xi \subseteq E^{\,\prime} \times E^{\,\prime}$, so as to form corresponding TTS-s $\sigma =
( E, E^{\, \prime}, T , \Xi )$. \\

{\bf Pre-Order on Topological Type Processes} \\

As mentioned, it is natural to assume the existence of a pre-order $\leq$ on the set
$E^{\,\prime}$ of topological type processes, or TTP-s, associated with the space $E$. Let us
therefore see how the above definitions leading to that in (TTS) are adapted in such a case. \\

First, the mapping $T : {\cal P} ( E ) \longrightarrow {\cal P} ( E^{\,\prime} )$ can
naturally be required to satisfy in addition to (TTP1) and (TTP2), also the condition \\

(TTP3)~~~ $ A \subseteq E,~ x^{\,\prime} \in T ( A ),~ y^{\,\prime} \in E^{\,\prime},~
        x^{\,\prime} \leq y^{\,\prime} ~~~\Longrightarrow~~~ y^{\,\prime} \in T ( A ) $ \\

Then the binary relation $\Xi \subseteq E^{\,\prime} \times E^{\,\prime}$ of Cauchy-ness, in
addition to (TTS1) - (TTS3), can naturally be required to satisfy the condition \\

(TTS4)~~~ $ x^{\,\prime} ~\Xi~ y^{\,\prime},~ x^{\,\prime} \leq u^{\,\prime},~
   y^{\,\prime} \leq v^{\,\prime} ~~~\Longrightarrow~~~ u^{\,\prime} ~\Xi~ v^{\,\prime} $ \\

for every $x^{\,\prime}, y^{\,\prime}, u^{\,\prime}, v^{\,\prime} \in E^{\,\prime}$. \\

Related or alternative such natural conditions, together with their consequences, are
presented in (3.1) - (3.3) and (3.20). \\

{\bf Topological Type Structures with Refinement, or TTSR-s} \\

In this way we can augment the definition of topological type structures in (TTS) with the
following definition of {\it topological type structures with refinement of their topological
type processes}, or in short TTSR-s, given by \\

(TTSR)~~~ $ \sigma ~=~ ( E, ( E^{\,\prime}, \leq ), T, \Xi ) $ \\

where $( E^{\,\prime}, \leq )$ is a pre-ordered set, and the conditions (TTP1) - (TTP3),
together with (TTS1) - (TTS4), are satisfied. \\

In view of (2.4), (2.5), it follows easily that we have \\

(2.9)~~~ $ x^{\,\prime} \in Cauchy ( \sigma ),~ y^{\,\prime} \in E^{\,\prime},~
 x^{\,\prime} \leq y^{\,\prime} ~~~\Longrightarrow~~~ y^{\,\prime} \in Cauchy ( \sigma ) $ \\

while for every $x \in E$, we have \\

(2.10)~~~ $ \begin{array}{l}
            x^{\,\prime} \in Conver ( \sigma, x ),~ y^{\,\prime} \in
                    E^{\,\prime},~ x^{\,\prime} \leq y^{\,\prime} ~~~\Longrightarrow \\ \\
   ~~~~~~~~~~~~~~~~~~~~~~~~~~~~~~~~~~~\Longrightarrow~~~ y^{\,\prime} \in Conver ( \sigma, x )
            \end{array} $ \\ \\

{\bf Two Basic Examples} \\

{\bf Usual Uniform Topologies as TTS-s.} Suppose we are given a usual uniform topology
$\upsilon$ on $E$. In this case we can, for instance, choose \\

(2.11)~~~ $ E^{\, \prime} = Fil ( E ) $ \\

Then the mapping (2.3) can be taken as \\

(2.12)~~~ $ T ( A ) ~=~ \{~ {\cal F} \in Fil ( E ) ~~|~~ A \in {\cal F} ~\},~~~
                                     \mbox{for}~~ A \subseteq E,~~ A \neq \phi $ \\

while $T ( \phi ) = \phi$, and we obtain the respective topological support on $E$ given by \\

(2.13)~~~ $ ( E, E^{\,\prime}, T ) ~=~ ( E, Fil ( E ), T ) $ \\

As for the binary relation (2.1), we can take it defined for any two filters ${\cal F, G}
\in Fil ( E )$, by \\

(2.14)~~~ $ {\cal F} ~\Xi_{\,\upsilon}~ {\cal G} ~~~\Longleftrightarrow~~~
               \left(~~ \begin{array}{l}
                   \exists~ x^{\#} \in E^{\#} : \\ \\
                   ~~~~ {\cal F} ~~\mbox{and}~~ {\cal G} ~~\mbox{converge in}~~ E^{\#}
                                           ~~\mbox{to}~~ x^{\#}
                \end{array} ~~\right) $ \\

where by $E^{\#}$ we denoted the usual completion of $E$ in the uniform topology $\upsilon$. \\

It is now easy to see that with the above choices \\

(2.15)~~~ $\sigma_\upsilon ~=~ ( E, E^{\, \prime}, T, \Xi_{\,\upsilon} ) $ \\

is a TTS on $E$, and $Cauchy~ ( \sigma_\upsilon )$ and $Conver~ ( \sigma_\upsilon, x )$, with
$x \in E$, are the respective sets of Cauchy and convergent filters on $E$ in the usual
uniform topology $\upsilon$ on $E$. \\
Also, $E$ is complete in the TTS $\sigma_\upsilon$, if and only if it is complete in the usual
uniform topology $\sigma$. \\

Furthermore, in this case the binary relation $\Xi_{\,\upsilon}$ is obviously {\it transitive},
which in the general definition of TTS-s is not required. \\

If we consider on $E^{\, \prime} = Fil ( E )$ the usual pre-order of filter refinement, which
in fact is a partial order, and it is simply given by the inclusion $\subseteq$ among filters,
then we have the stronger result \\

(2.15$^{\,\prime}$)~~~ $\sigma_\upsilon ~=~
                        ( E, ( E^{\, \prime}, \subseteq ), T, \Xi_{\,\upsilon} ) $ \\

is a TTSR on $E$. \\

{\bf Usual Topologies as TTS-s.} Now, let us be given a usual topology $\tau$ on $E$. Then we
can take on $E$ the same topological support as in (2.13). For the binary relation (2.1), we
take \\

(2.16)~~~ $ \Xi_{\,\tau} ~=~ \left \{~ ( {\cal F, G } ) \in Fil ( E) \times Fil ( E ) ~~
                 \begin{array}{|l}
                    \exists~ x \in E : \\ \\
                    ~~~~ {\cal F, G} ~~\stackrel{\tau} \longrightarrow~~ x
                 \end{array} ~ \right \} $ \\

where ${\cal F} ~~\stackrel{\tau} \longrightarrow~~ x$ means that the filter ${\cal F}$ on $E$
converges to $x$ in the usual topology $\tau$ on $E$. \\

It follows easily that with the above choices \\

(2.17)~~~ $ \sigma_\tau ~=~ ( E, E^{\, \prime}, T, \Xi_{\,\tau} ) $ \\

is a TTS on $E$, with $Conver~ ( \sigma_\tau, x )$ being the set of filters on $E$ convergent
to $x$ in the usual topology $\tau$ on $E$. \\
Clearly, we also have \\

(2.18)~~~ $ Cauchy~ ( \sigma_\tau ) ~=~ \bigcup_{x \in E}~ Conver~ ( \sigma_\tau, x ) $ \\

which means that $\sigma_\tau$ is a complete TTS on $E$. Therefore this particular way, namely,
$\tau \longmapsto \sigma_\tau$, to associate a TTS on $E$ to a usual topology $\tau$ on $E$ is
rather trivial. \\

We also note that the binary relation $\Xi_{\,\tau}$ is {\it transitive}. \\

If again, we consider on $E^{\, \prime} = Fil ( E )$ the usual partial order of filter
refinement, which in fact is simply given by the inclusion $\subseteq$ among filters, then we
have the stronger result \\

(2.17$^{\,\prime}$)~~~ $ \sigma_\tau ~=~
                        ( E, ( E^{\, \prime}, \subseteq ), T, \Xi_{\,\tau} ) $ \\

is a TTSR on $E$. \\

{\bf Associations of Alternative TTS-s.} Obviously, to a given usual topology $\tau$ on $E$
one may as well associate other TTS-s than the above $\sigma_\tau$. And this may, among others,
be done in order to obtain a completion of $E^{\#}$ which, unlike as in (2.18), may be larger
than $E$ itself. One way to do that is by considering on $E$ certain TTS-s \\

(2.19)~~~ $ \sigma ~=~ ( E, E^{\,\prime}, T, \Xi ) $ \\

for which we have \\

(2.20)~~~ $ \Xi_{\,\tau} ~\subseteq~ \Xi $ \\

In such a case obviously we will have \\

(2.21)~~~ $ \begin{array}{l}
                  Cauchy~ ( \sigma_\tau ) ~\subseteq~ Cauchy~ ( \sigma ) \\ \\
                  Conver~ ( \sigma_\tau, x ) ~\subseteq~ Conver~ ( \sigma_\tau, x ),
                                       ~~\mbox{for}~~ x \in E
             \end{array} $ \\

while at the same time it may happen that $\sigma$ itself is complete, that is, it satisfies
(2.7). \\

{\bf The Nonrigidity of TTSR-s} \\

It is important to note that the TTSR-s defined above are {\it nonrigid} mathematical
structures. Indeed, given $\sigma = ( E, ( E^{\,\prime}, \leq ), T, \Xi )$ a TTSR on a set $E$,
it involves {\it four} defining entities, namely, a set $E^{\,\prime}$, a pre-order $\leq$ on
that set, a mapping $T$, and a binary relation $\Xi$ on $E^{\,\prime}$. And as shown by the
respective conditions required about these four entities, those conditions are merely {\it
compatibility} conditions, thus they are {\it not} determining any one of these four entities
in function of the other three. \\
This is in sharp contradistinction to the case of the usual HKB concept of topology, where for
instance, the open sets are supposed to determine uniquely all the other topological
entities. \\

Therefore, as mentioned earlier, it is precisely this {\it nonrigid} aspect of the TTSR-s
which offers the possibility of a suitable extension of the usual HKB concept of
topology. \\ \\

{\large \bf 3. Topological Type Structures with Pre-ordered \\
\hspace*{0.55cm} Topological Type Processes Can Lead Back to \\
\hspace*{0.55cm} Usual HKB Topologies} \\

Here we shall show that, under rather general and natural conditions, TTS-s which have
pre-ordered TTP-s can {\it lead back} to usual HKB topologies. And in fact, we shall give two
such constructions, Rosinger [7, pp. 141-151]. Needless to say, there may be other similar
constructions as well. \\

These constructions give an indication of the fact that the concept of TTS as defined in this
work and originated in Rosinger [7], is {\it not} unduly abstract or general. \\
Also, they show how natural is to consider a pre-order on TTP-s. \\
On the other hand, one should note that the mentioned natural constructions leading from
TTS-s back to usual HKB topologies do {\it not} necessarily mean that the former are reduced
to the latter, since as mentioned, the former will typically be more general due, among others,
to their {\it nonrigid} structure. Moreover, one can only do so under certain conditions which,
albeit, appear to be quite natural. \\
This argument is further supported by the fact that, as seen next, from a given TTS one can
go back in more than one natural way to usual HKB topologies. \\

{\bf Compatible Pre-ordered Topological Type Processes} \\

In order to illustrate the conceptual nonrigidity involved in the present approach to
pseudo-topologies, here we deal with the possible connection between TTS-s and pre-orders on
TTP-s in an alternative manner to that leading to the definition of TTSR-s given in (TTSR)
above. \\

It is important to note that usual topological type processes, or TTP-s, have a {\it two fold}
connection with pre-orders, Rosinger [7, pp. 141-151], namely :

\begin{itemize}

\item each such TTP on a given space $E$ is a mapping of a pre-ordered directed set into a
certain suitably given set,

\item there is a further pre-order on the set of all such TTP-s associated to the space $E$.

\end{itemize}

As mentioned earlier, certainly the above is the case when, for instance, TTP-s are given by
nets or filters. Indeed, nets on a given nonvoid set $E$ are mappings of directed pre-orders
into $E$, while two such nets can be compared with one another when one is a subnet of the
other. Similarly, filters ${\cal F}$ on $E$ can be seen as mappings ${\cal F} \ni F
\longmapsto F \in {\cal P} ( E )$. Thus they are nets with values in ${\cal P} ( E )$, since
${\cal F}$ is directed with respect to the partial order $\supseteq$. Furthermore, two such
filters on $E$ can be compared with one another when one is more fine than the other. \\

Let $\sigma = ( E, E^{\, \prime}, T , \Xi )$ be any TTS on the nonvoid set $E$ and let $\leq$
be any pre-order on $E^{\, \prime}$, which we recall, is the set of TTP-s of $\sigma$. \\

We call $\sigma = ( E, E^{\, \prime}, T , \Xi )$ and $( E^{\, \prime}, \leq )$ {\it
compatible}, if and only if the following three conditions hold \\

(3.1)~~~ $ \begin{array}{l}
               \forall~~~ A \subseteq E,~ x^{\, \prime} \in T ( A ),~
                                             y^{\, \prime} \in E^{\, \prime} ~: \\ \\
               ~~~~~ x^{\, \prime} ~\leq~ y^{\, \prime} ~~~\Longrightarrow~~~
                                    y^{\, \prime} \in T ( A )
           \end{array} $ \\ \\ \\

(3.2)~~~ $ \begin{array}{l}
               \forall~~~ x \in E,~ x^{\, \prime} \in Conver~ ( \sigma, x ),~
                                           y^{\, \prime} \in E^{\, \prime} ~: \\ \\
               ~~~~~ x^{\, \prime} ~\leq~ y^{\, \prime} ~~~\Longrightarrow~~~
                                    y^{\, \prime} \in Conver~ ( \sigma, x )
           \end{array} $ \\ \\ \\

(3.3)~~~ $ \begin{array}{l}
                \forall~~~ x^{\, \prime} \in Cauchy~ ( \sigma ),~
                                           y^{\, \prime} \in E^{\, \prime} ~: \\ \\
                ~~~~~ x^{\, \prime} ~\leq~ y^{\, \prime} ~~~\Longrightarrow~~~
                  \left ( ~~ \begin{array}{l}
                              y^{\, \prime} \in Cauchy~ ( \sigma ) \\ \\
                              ( x^{\, \prime}, y^{\, \prime} ) \in~ \Xi
                           \end{array} ~~ \right )
            \end{array} $ \\ \\

It is easy to see that in case the TTP-s are given by nets or filters, then the above
compatibility conditions (3.1) - (3.3) are indeed satisfied. \\

Let us further note that (3.1) is but an equivalent reformulation of the earlier (TTP3), while
in view of (TTS3), the condition (3.3) is implied by (TTS4). \\

{\bf Closed and Open Subsets, Neighbourhoods} \\

Let us for a moment return to general TTS-s $\sigma ~=~ ( E, E^{\, \prime}, T , \Xi )$. A
subset $A \subseteq E$ is called $\sigma$-{\it closed}, if and only if \\

(3.4)~~~ $ \begin{array}{l}
               \forall~~~ x^{\, \prime} \in T ( A ),~ x \in E ~: \\ \\
               ~~~~~ x^{\, \prime} \in Conver~ ( \sigma, x ) ~~~\Longrightarrow~~~ x \in A
           \end{array} $ \\ \\

We denote by $Cl~ ( \sigma )$ the set of all $\sigma$-closed subsets of $E$. \\

A subset $A \subseteq E$ is called $\sigma$-{\it open}, if and only if $E \setminus A$ is
$\sigma$-closed. And we denote by $Op~ ( \sigma )$ the set of all $\sigma$-open subsets of
$E$. \\

It is easy to see that \\

(3.5)~~~ \mbox{arbitrary intersections of}~ $\sigma$-closed
                                  ~\mbox{sets are}~ $\sigma$-closed \\

(3.6)~~~ \mbox{arbitrary unions of}~ $\sigma$-open
                                  ~\mbox{sets are}~ $\sigma$-open \\

(3.7)~~~ $\phi$ ~\mbox{and}~ $E$ ~\mbox{are both}~
                                 $\sigma$-closed ~\mbox{and}~ $\sigma$-open \\

Consequently, we can define the operation of $\sigma$-{\it closure} as follows \\

(3.8)~~~ $ E \supseteq A ~~~\longmapsto~~~ cl_\sigma A ~=~
            \bigcap_{\, B\, \in\, Cl\, ( \sigma ),~ B \supseteq A}~ B \in Cl~ ( \sigma ) $ \\

and of $\sigma$-{\it interior} by \\

(3.9)~~~ $ E \supseteq A ~~~\longmapsto~~~ in_\sigma A ~=~
            \bigcup_{\, B\, \in\, Op\, ( \sigma ),~ B \subseteq A}~ B \in Op~ ( \sigma ) $ \\

It is easy to see that the following properties hold for all subsets $A \subseteq E$ \\

(3.10)~~~ $ A \in Cl~ ( \sigma ) ~~~\Longleftrightarrow~~~ A ~=~ cl_\sigma A $ \\

(3.11)~~~ $ A ~\subseteq~ cl_\sigma A ~=~ cl_\sigma~ cl_\sigma A $ \\

(3.12)~~~ $ A \in Op~ ( \sigma ) ~~~\Longleftrightarrow~~~ A ~=~ op_\sigma A $ \\

(3.13)~~~ $ A ~\supseteq~ op_\sigma A ~=~ op_\sigma~ op_\sigma A $ \\

Also, for $A \subseteq B \subseteq E$ we have \\

(3.14)~~~ $ cl_\sigma A ~\subseteq~ cl_\sigma B,~~~ op_\sigma A ~\subseteq~ op_\sigma B $ \\

while \\

(3.15)~~~ $ cl_\sigma \phi ~=~ op_\sigma \phi ~=~ \phi,~~~
                                 cl_\sigma E ~=~ op_\sigma E ~=~ E $ \\

Let us now define the set of $\sigma$-{\it neighbourhoods} of a point $x \in E$ by \\

(3.16)~~~ $ {\cal V}_\sigma ( x ) ~=~
              \{~ A \subseteq E ~~|~~ \exists~~ B \in Op~ ( \sigma ) ~:~
                                                         x \in B \subseteq A ~\} $ \\

It follows easily that for a subset $A \subseteq E$, we have \\

(3.17)~~~ $ A \in Op~ ( \sigma ) ~~~\Longleftrightarrow~~~
                     \forall~~ x \in A ~:~ A \in {\cal V}_\sigma ( x ) $ \\

{\bf From TTS-s Back to HKB Topologies} \\

As seen above in (3.5) - (3.17), within arbitrary TTS-s $\sigma = ( E, E^{\, \prime}, T ,
\Xi )$ one can naturally define corresponding concepts of closed and open sets, closure and
interior operators, as well as neighbourhoods of points in $E$. And in view of (3.5) and (3.6),
the only thing which prevents such associated concepts from being identical with those in
usual HKB topologies is that, in general, finite unions of closed sets need not be closed, and
thus, finite intersections of open sets need not be open. \\

Here we show that in case TTS-s $\sigma ~=~ ( E, E^{\, \prime}, T , \Xi )$ are {\it compatible}
with suitably associated pre-orders $( E^{\, \prime}, \leq )$, then the above association in
(3.5) - (3.17) with corresponding topological concepts does in fact lead to usual HKB
topologies. \\

For that purpose, it is useful to introduce the following $\sigma$-{\it adherence} operator,
Rosinger [7, pp. 141-151] \\

(3.18)~~~ $ T_\sigma : {\cal P} ( E ) \times {\cal P} ( E )
                          ~~~\longrightarrow~~~ {\cal P} ( E^{\, \prime} ) $ \\

defined for $A, B \subseteq E$ by \\

(3.19)~~~ $ T_\sigma ( A, B ) ~=~ \{~ x^{\, \prime} \in T ( B ) ~~|~~
            \exists~~~ y^{\, \prime} \in T ( A ) ~:~ x^{\, \prime} \leq y^{\, \prime} ~\} $ \\

In terms of this operator, the earlier compatibility conditions (3.1) - (3.3) are now
strengthened as follows. \\

We call $\sigma ~=~ ( E, E^{\, \prime}, T , \Xi )$ and $( E^{\, \prime}, \leq )$ {\it strongly
compatible}, if and only if in addition to (3.1) - (3.3), the following condition holds : \\

(3.20)~~~ $ \begin{array}{l}
                \forall~~~ A, B \subseteq E ~: \\ \\
                ~~~~~ T ( A \cup B ) ~\subseteq~
                                       T_\sigma ( A, A \cup B ) \cup T_\sigma ( B, A \cup B )
             \end{array} $ \\

Once again it is easy to see that in case the TTP-s are given by nets or filters, then the
above compatibility condition (3.20) is indeed satisfied. \\

Now we can obtain the following \\

{\bf Proposition} ( Rosinger [7, pp. 145,146] )\\

Given $\sigma$ and $( E^{\, \prime}, \leq )$ strongly compatible, then \\

(3.21)~~~ finite unions of $\sigma$-closed sets are $\sigma$-closed \\

(3.22)~~~ finite intersections of $\sigma$-open sets are $\sigma$-open \\

Also, for every $x \in E,~ x^{\, \prime} \in E^{\, \prime}$, we have \\

(3.23)~~~ $ x^{\, \prime} \in Conver ( \sigma, x ) \Longrightarrow
                   \left (  \begin{array}{l}
                                \forall~ A \in {\cal V}_\sigma ( x ) ~: \\ \\
                                \exists~ y^{\, \prime} \in Conver ( \sigma, x )
                                                   \cap T ( A ) ~: \\ \\
                                ~~~ x^{\, \prime} \leq y^{\, \prime}
                              \end{array}  \right ) $ \\

{\bf Proof} \\

For (3.21), let $A, B$ be $\sigma$-closed in $E$, $x^{\, \prime} \in T ( A \cup B )$ and $x
\in E$, such that $x^{\, \prime} \in Conver ( \sigma, x )$. Then we show that $x \in A \cup B$.
Thus in view of (3.4), that will result in $A \cup B$ being $\sigma$-closed. \\

Indeed $x^{\, \prime} \in T ( A \cup B )$ and (3.20) give $x^{\, \prime} \in T_\sigma ( A, A
\cup B )$, or $x^{\, \prime} \in T_\sigma ( B, A \cup B )$. In case the first relation holds,
then according to (3.19) it follows that $~\exists~~ y^{\, \prime} \in T ( A ) ~:~ x^{\,
\prime} \leq y^{\, \prime}$. Hence $x^{\, \prime} \in Conver ( \sigma, x )$ and (3.2) give
$y^{\, \prime} \in Conver ( \sigma, x )$. But then $y^{\, \prime} \in T ( A )$ and $A$ being
$\sigma$-closed imply $x \in A$, and thus indeed (3.21) holds. Consequently, so does (3.22). \\

We show now (3.23). Let $A \in {\cal V}_\sigma ( x )$, then (3.16) gives $B \in Op~ ( \sigma
)$, with $x \in B \subseteq A$. But obviously $x^{\, \prime} \in Conver ( \sigma, x )
\subseteq E^{\, \prime} = T ( E ) \subseteq T_\sigma ( B, E ) \cup T_\sigma ( E \setminus B,
E )$. Assume now that $x^{\, \prime} \in T_\sigma ( E \setminus B, E )$. And then $x^{\,
\prime} \in Conver ( \sigma, x )$ and $E \setminus B$ being $\sigma$-closed imply $x \in E
\setminus B$, which in view of the above is absurd. \\
Thus $x^{\, \prime} \in T_\sigma ( B, E )$ which according to (3.19) gives $y^{\, \prime} \in
T ( B )$, with $x^{\, \prime} \leq y^{\, \prime}$. But then (3.2) yields $y^{\, \prime} \in
Conver ( \sigma, x )$, and the proof of (3.23) is completed since $y^{\, \prime} \in T ( B )
\subseteq T ( A )$, as we have $A \subseteq B$.

\hfill $\Box$ \\

In view of (3.21), (3.22), it follows that under the conditions of the above Proposition, we
can associate with the TTS $\sigma ~=~ ( E, E^{\, \prime}, T , \Xi )$ a usual HKB topology
$\tau_\sigma$ on $E$ which will have the {\it same} closed and open sets. \\

Here, however, one should note that the associated topology $\tau_\sigma$ does {\it not} in
general express all the structure or information in the initial TTS $\sigma ~=~
( E, E^{\, \prime}, T , \Xi )$. For instance, it need not fully express the {\it uniform}
aspects involved in $\Xi$. \\

{\bf Associating Compatible Topological Type Processes} \\

In view of the above, there is an interest in associating to {\it arbitrary} TTS-s $\sigma ~=~
( E, E^{\, \prime}, T , \Xi )$ certain {\it compatible} pre-ordered TTP-s, in order to be
further able to associate usual HKB topologies, in case the respective pre-ordered TTP-s turn
out to be strongly compatible as well. \\
Here we shall give {\it two} such constructions in which we associate compatible pre-ordered
TTP-s with {\it arbitrary} TTS-s. \\
As for conditions when such associated pre-ordered TTP-s are in fact strongly compatible with
the given TTS-s, this is an issue which can be dealt with separately. \\

{\bf The First Association of Pre-ordered TTP-s} \\

Given any $\sigma ~=~ ( E, E^{\, \prime}, T , \Xi )$ which is a TTS on the nonvoid set $E$. We
associate with $\sigma$ the following TTS on the same set $E$, as defined by \\

(3.24)~~~ $ \sigma_P ~=~ ( E, E^{\, \prime}_P, T_P , \Xi_P ) $ \\

Here \\

(3.25)~~~ $ E^{\, \prime}_P ~=~
                   {\cal P} ( E^{\, \prime} ) $ \\

while the pre-order, in fact, a partial order $\leq_P$ on $ E^{\, \prime}_P$ is defined by \\

(3.26)~~~ $ x^{\, \prime}_P ~\leq_P~ y^{\, \prime}_P
                   ~~~\Longleftrightarrow~~~ x^{\, \prime}_P ~\supseteq~ y^{\, \prime}_P $ \\

for every $x^{\, \prime}_P, y^{\, \prime}_P \in E^{\, \prime}_P$. \\

Further, the mapping \\

(3.27)~~~ $ T_P : {\cal P} ( E ) ~~\longrightarrow~~ {\cal P} (  E^{\, \prime}_P ) $ \\

is defined for $A \subseteq E$, by \\

(3.28)~~~ $ T_P ( A ) ~=~ {\cal P} ( T ( A ) ) $ \\

Finally, $\Xi_P$ is defined as being the set of all pairs $( x^{\, \prime}_P, y^{\, \prime}_P ) \in
E^{\, \prime}_P \times E^{\, \prime}_P$ such that \\

(3.29)~~~ $ \begin{array}{l}
                \forall~~~ x^{\, \prime}, y^{\, \prime} \in
                                 x^{\, \prime}_P \cup y^{\, \prime}_P ~: \\ \\
                ~~~~~ ( x^{\, \prime}, y^{\, \prime} ) \in \Xi
            \end{array} $ \\ \\

The result of interest for us is in the following easy to prove \\

{\bf Proposition} ( Rosinger [7, p. 148] ) \\

Given any $\sigma ~=~ ( E, E^{\, \prime}, T , \Xi )$ which is a TTS on the nonvoid set $E$.
Then $\sigma_P$ and $( E^{\, \prime}, \leq_P )$ are compatible. \\

{\bf The Second Association of Pre-ordered TTP-s} \\

Let $\eta = ( E, E^{\, \prime}, T )$ be a topological support on the nonvoid set $E$, see
(2.8). We define the mapping \\

(3.30)~~~ $ E^{\, \prime} \ni x^{\, \prime}
               ~\longmapsto~ {\cal W}_\eta ( x^{\, \prime} ) ~=~
                   \{~ \phi \neq A ~\subseteq~ E ~~|~~ x^{\, \prime} \in T ( A ) ~\} $ \\

The we have for $x^{\, \prime} \in E^{\, \prime}$ the relations \\

(3.31)~~~ $ \phi \notin {\cal W}_\eta ( x^{\, \prime} ) ~\neq~ \phi $ \\

while for $\phi \neq A \subseteq B \subseteq E$, we have \\

(3.32)~~~ $ A \in {\cal W}_\eta ( x^{\, \prime} ) ~~\Longrightarrow~~
                                    B \in {\cal W}_\eta ( x^{\, \prime} ) $ \\

Consequently, ${\cal W}_\eta ( x^{\, \prime} )$, with $x^{\, \prime} \in E^{\, \prime}$, are
filters on $E$, if and only if \\

(3.33)~~~ $ T ( A ) \cap T ( B ) ~\subseteq~ T ( A \cap B ),~~
                                \mbox{with}~ \phi \neq A, B \subseteq E $ \\

Finally, we define the partial order $\leq_\eta$ on $E^{\, \prime}$, by \\

(3.34)~~~ $ x^{\, \prime} ~\leq_\eta~ y^{\, \prime} ~~~\Longleftrightarrow~~~
             {\cal W}_\eta ( x^{\, \prime} ) ~\subseteq~ {\cal W}_\eta ( y^{\, \prime} ) $ \\

where $x^{\, \prime}, y^{\, \prime} \in E^{\, \prime}$. \\

Now we obtain the following \\

{\bf Proposition} ( Rosinger [7, p. 151] ) \\

Given $\sigma ~=~ ( E, E^{\, \prime}, T , \Xi )$ a TTS on the nonvoid set $E$, and let  $\eta
= ( E, E^{\, \prime}, T )$ be its topological support. \\

Then $\sigma$ and $( E^{\, \prime}, \leq_\eta )$ are compatible, if and only if (3.2) and (3.3)
hold. \\

{\bf Proof} \\

It suffices to show that (3.1) is satisfied. This however follows immediately from
(3.34). \\ \\

{\large \bf 4. The Particular Case of Beattie-Butzmann \\
\hspace*{0.55cm} Convergence and Uniform Convergence Spaces} \\

{\bf Convergence Spaces} \\

As mentioned, one of the most successful and seminal ways to go beyond usual topological or
uniform spaces was developed in a series of publications by Beattie  \& Butzmann. Let us for
convenience recall here their first basic definition, Beattie  \& Butzmann [p. 2]. \\

{\bf Definition.} On a given set $E$ we call {\it convergence structure} any mapping \\

(4.1)~~~ $ \lambda : E ~\longrightarrow~ {\cal P} ( Fil ( E ) ) $ \\

of $E$ into the powers set of all filters on $E$, which satisfies the following three
conditions. First, we have \\

(4.2)~~~ $ {\cal U}_{\,x} \in \lambda ( x ), ~~\mbox{for}~~ x \in E $ \\

where ${\cal U}_{\,x}$ denotes the filter, in fact, an ultrafilter, generated on $E$ by the
one point set $\{ x \}$. Second, we have for all $x \in E$ \\

(4.3)~~~ $ {\cal F, G} \in \lambda ( x ) ~~~\Longrightarrow~~~
                             {\cal F} \cap {\cal G} \in \lambda ( x ) $ \\

Finally, for all $x \in E$, we also have \\

(4.4)~~~ $ {\cal F} \in \lambda ( x ),~ {\cal G} \in Fil ( E ),~ {\cal G}
           ~\supseteq~ {\cal F} ~~~\Longrightarrow~~~ {\cal G} \in \lambda ( x ) $ \\

In such a case one calls $( E, \lambda )$ a {\it convergence space} on $E$. Furthermore,
for every $x \in E$ and ${\cal F} \in \lambda ( x )$, we say that the filter ${\cal F}$ {\it
converges to} $x$ in the convergence structure $( E, \lambda )$, and thus write ${\cal F}
\stackrel{\lambda} \longrightarrow x$. \\

{\bf Motivation.} The motivation of the above definition is rather straightforward. Let $\tau$
be any usual topology on $E$. Then the mapping $\lambda$ in (4.1) is defined by \\

(4.5)~~~ $ \lambda ( x ) ~=~ \{~ {\cal F} \in Fil ( E ) ~~|~~
                                      {\cal V}_x \subseteq {\cal F} ~\},~~ x \in E $ \\

where ${\cal V}_x$ denotes the filter of {\it neighbourhoods} of $x$ in the topology $\tau$. \\

{\bf Inclusion Among TTS-s.} We show now the simple way in which the convergence spaces $( E,
\lambda )$ are {\it particular} cases of the TTS-s defined in section 1. \\
Indeed, there is a natural way one can associate with each convergence space $( E, \lambda )$
a TTS~ $\sigma_{\,\lambda}$ on $E$. Namely, we take \\

(4.6)~~~ $ \sigma_{\,\lambda} ~=~ ( E, E^{\,\prime}, T, \Xi_{\,\lambda} ) $ \\

where the topological support $( E, E^{\,\prime}, T)$ is the same with the one in (2.13),
while similar to (2.16), we take \\

(4.7)~~~ $ \Xi_{\,\lambda} ~=~ \left \{~ ( {\cal F, G } ) \in Fil ( E) \times Fil ( E ) ~~
                 \begin{array}{|l}
                    ~~\exists~ x \in E : \\ \\
                    ~~~~~~ {\cal F, G} ~~\stackrel{\lambda} \longrightarrow~~ x
                 \end{array} ~ \right \} $ \\

It is easy to see that $\sigma_{\,\lambda} ~=~ ( E, E^{\,\prime}, T, \Xi_{\,\lambda} )$ is
indeed a TTS on $E$. \\

Furthermore, as in (2.17$^{\,\prime}$) we can associate with the convergence space $( E,
\lambda )$ not only a TTS, but also a TTSR. \\

{\bf Uniform Convergence Spaces} \\

The second basic definition in Beattie \& Butzmann, see [p. 61], is the following \\

{\bf Definition.} On a given set $E$ a {\it uniform convergence structure} is any set {\large
\bf U} of filters on $E \times E$, that is, {\large \bf U} $\subseteq Fil ( E \times E )$,
such that the following five conditions hold \\

(4.8)~~~ ${\cal U}_{( x, x )} \in$~{\large \bf U},~~ for $x \in E$ \\

(4.9)~~~ ${\cal U,~ V} \in$~{\large \bf U} $~~~\Longrightarrow~~~
                         {\cal U} \bigcap {\cal V} \in$~{\large \bf U} \\

(4.10)~~~ ${\cal U} \in$~{\large \bf U},~ ${\cal V} \in Fil ( E \times E ),~
  {\cal V} \supseteq {\cal U} ~~~\Longrightarrow~~ {\cal V} \in$~ {\large \bf U} \\

(4.11)~~~ ${\cal U} \in$~{\large \bf U}
                       $~~~\Longrightarrow~~~ {\cal U}^{-1} \in$~{\large \bf U} \\

and \\

(4.12)~~~  ${\cal U,~ V} \in$~{\large \bf U},~~ ${\cal U} \circ {\cal V} \in
                   Fil ( E \times E ) ~~~\Longrightarrow~~~
                       {\cal U} \circ {\cal V} \in$~{\large \bf U} \\

Above, as usual, we denoted \\

~~~~~~ $ {\cal U}_{( x, x )} = \{~ C \subseteq E \times E ~|~ ( x, x ) \in C ~\} $ \\

~~~~~~ $ {\cal U}^{-1} = \{~ C^{-1} ~|~ C \in {\cal U} ~\} $ \\

where \\

~~~~~~ $ C^{-1} = \{~ ( y, x ) ~|~ ( x, y ) \in C ~\} $ \\

and lastly \\

~~~~~~ $ {\cal U} \circ {\cal V} = \{~ C \subseteq E \times E ~|~ \exists~ A \in {\cal U},~
           B \in {\cal V} : A \circ B \subseteq C ~\} $ \\

where \\

~~~~~~ $ A \circ B = \{~ ( x, z ) ~|~
                          \exists~ y \in E : ( x, y ) \in A,~ ( y, z ) \in B ~\} $ \\

In this case one calls ( $E$, {\large \bf U} ) a {\it uniform convergence space}. \\

{\bf Inclusion Among TTS-s.} We can now show the simple way in which the uniform convergence
spaces ( $E$, {\large \bf U} ) are again {\it particular} cases of the TTS-s defined in
section 1. \\
There is, indeed, a natural way one can associate with each uniform convergence space ( $E$,
{\large \bf U} ) a TTS~ $\sigma_{\,{\large \bf U}}$ on $E$. Namely, we take \\

(4.13)~~~ $ \sigma_{\,{\large \bf U}} ~=~ ( E, E^{\,\prime}, T, \Xi_{\,{\large \bf U}} ) $ \\

where the topological support $( E, E^{\,\prime}, T)$ is the same with the one in (2.13),
while \\

(4.14)~~~ $ \Xi_{\,{\large \bf U}} ~=~ \{~ ( {\cal F, G } ) \in Fil ( E) \times Fil ( E )
               ~~|~~ {\cal F} \times {\cal G} \in$~ {\large \bf U} ~\} \\

where as usual, we denoted \\

~~~~~~ $ {\cal F} \times {\cal G} ~=~ \{~ C \subseteq E \times E ~~|~~
             \exists~ A \in {\cal F},~ B \in {\cal G} : A \times B \subseteq C ~\} $ \\

It is easy to see that $\sigma_{\,{\large \bf U}} ~=~ ( E, E^{\,\prime}, T,
\Xi_{\,{\large \bf U}} )$ is indeed a TTS on $E$. \\


\end{document}